\documentclass[letterpaper, 10 pt, conference]{ieeeconf}  

\usepackage[ansinew]{inputenc}
\usepackage{float}
\usepackage{graphicx}
\usepackage{fancyhdr}
\usepackage{amsmath,amssymb}
\usepackage{color}
\usepackage{verbatim}
\usepackage{scrtime}
\usepackage{float}

\usepackage{subfigure}
\usepackage{url}
\usepackage{multirow}
\usepackage[sort]{cite}

\usepackage[usenames,dvipsnames]{xcolor}
\usepackage{hyperref}
\hypersetup{
     colorlinks  = true,
     citecolor  = blue,
     linkcolor=blue
}

\usepackage{algpseudocode}
\usepackage{algorithm}
\algtext*{EndIf}
\algnewcommand\algorithmicinput{\textbf{Input:}}
\algnewcommand\Input{\item[\algorithmicinput]}
\algnewcommand\algorithmicoutput{\textbf{Ouput:}}
\algnewcommand\Output{\item[\algorithmicoutput]}

\newcommand{\change}[1]{#1}

\newcommand{\trace}{{\rm Trace}}
\newcommand{\Stiefel}[2]{{\mathrm{St}({#1},{#2})}}
\newcommand{\Gr}[2]{{\mathrm{Gr}({#1},{#2})}}
\newcommand{\OG}[1]{{\mathcal{O}({#1})}}
\newcommand{\Sym}{{\mathrm{Sym}}}
\newcommand{\Skew}{{\mathrm{Skew}}}

\newcommand{\argmin}{\operatornamewithlimits{arg\,min}}

\newcommand{\grad}{\mathrm{grad}}

\newcommand{\mat}[1]{{\bf #1}}
\newcommand{\subject}{\mathrm{subject\  to}}

\newcommand{\rtwomc}{R$2$MC}
\newcommand{\rmc}{R$3$MC}

\IEEEoverridecommandlockouts                              
\title{\LARGE \bf
{\rmc}: A Riemannian three-factor algorithm for \\low-rank matrix completion$^\ast$\thanks{$^\ast$This paper presents research results of the Belgian Network DYSCO (Dynamical Systems, Control, and Optimization), funded by the Interuniversity Attraction Poles Programme, initiated by the Belgian State, Science Policy Office. The scientific responsibility rests with its authors. This work was also partly supported by the Belgian FRFC (Fonds de la Recherche Fondamentale Collective). Bamdev Mishra is a research fellow of the Belgian National Fund for Scientific Research (FNRS).}}

\author{Bamdev Mishra$^\dagger$\thanks{$^\dagger$ Department of Electrical Engineering and Computer Science, University of Li\`ege, 4000 Li\`ege,
Belgium (\{B.Mishra, R.Sepulchre\}@ulg.ac.be).}
 and Rodolphe Sepulchre$^{\S \dagger}$\thanks{$^\S$University of Cambridge, Department of Engineering, Trumpington Street, Cambridge CB2 1PZ, UK
(R.Sepulchre@eng.cam.ac.uk).}
}

\begin{document}

\maketitle
\thispagestyle{plain}
\pagestyle{plain}

\begin{abstract}
We exploit the versatile framework of Riemannian optimization on quotient manifolds to develop {\rmc}, a nonlinear conjugate-gradient method for low-rank matrix completion. The underlying search space of fixed-rank matrices is endowed with a novel Riemannian metric that is tailored to the least-squares cost. Numerical comparisons suggest that {\rmc} robustly outperforms state-of-the-art algorithms across different problem instances, especially \change{those that} combine scarcely sampled and ill-conditioned data.
\end{abstract}

\section{Introduction}
We address the problem of low-rank matrix completion when the rank is a priori known or estimated. Given a rank-$r$ matrix $\mat{X}^{\star}\in\mathbb{R}^{n\times m}$ whose entries $\mat{X}^{\star}_{ij}$ are only known for some indices $(i,j)\in\Omega$, where $\Omega$ is a subset of the complete set of indices $\{(i,j):i\in\{1,\dots,n\}\text{ and }j\in\{1,\dots,m\}\}$, the \emph{fixed-rank matrix completion problem} for rank $r$ is formulated as
\begin{equation}\label{eq:matrix_completion}
\begin{array}{llll}
	\min\limits_{\mat{X} \in \mathbb{R}^{n\times n}}
		&	\frac{1}{|\Omega|}\|\mathcal{P}_{\Omega}(\mat{X}) - \mathcal{P}_{\Omega}(\mat{X}^{\star})\|_F^2 \\
	 \subject &  {\rm rank}(\mat X) = r,
\end{array}
\end{equation}
where the operator $\mathcal{P}_{\Omega}(\mat{X})_{ij}=\mat{X}_{ij}$ if $(i,j) \in \Omega$ and $\mathcal{P}_{\Omega}(\mat{X})_{ij}=0$ otherwise,  $|\Omega|$ is the cardinality of the $\Omega$ (equal to the number of known entries), and the norm $\|\cdot \|_F$ is the Frobenius norm. Problem (\ref{eq:matrix_completion}) amounts to minimizing the data fitting error with the known entries in $\mat{X}^\star$ while constraining the rank to $r$. \change{Not surprisingly, Problem (\ref{eq:matrix_completion}) and its (many) variants find applications} in system identification \cite{markovsky08a}, control system \cite{benner13a}, computer vision \cite{joulin10a}, machine learning \cite{rennie05a, meyer11a}, to name a few.

The case of interest is when $r\ll \min(m, n)$, i.e., the rank is much smaller than the matrix dimensions. Solvability of (\ref{eq:matrix_completion}) under the low-rank assumption has been studied in \cite{candes09b, keshavan10a}, while a number of computationally efficient algorithms have been proposed in \cite{rennie05a, keshavan10a,  cai10a, jain10a, meyer11a,  boumal11a, mishra13a, mishra12a, mishra14a, ngo12a, vandereycken13a}, among others.

The problem (\ref{eq:matrix_completion}) has two fundamental structures. First, the \emph{least-squares} structure of the cost function. Second, the fixed-rank constraint that has the structure of a \change{\emph{quotient matrix manifold}} \cite{mishra14a}. We denote the set of $n\times m$ matrices of rank $r$ by $\mathbb{R}_r^{n \times m}$. A popular way to characterize fixed-rank matrices is through fixed-rank \emph{matrix factorizations} \cite{mishra14a}. However, most matrix factorizations have invariance properties that make the factorizations \emph{non-unique}. And in many cases the set $\mathbb{R}_r^{n\times m}$ is identified with structured differentiable quotient manifolds \cite{mishra14a}. Our focus in \change{the present} paper is on exploiting these particular structures in order to develop an efficient algorithm that \change{scales} to large-scale dimensional data by using manifold optimization techniques \cite{absil08a}.

The manifold optimization framework of \cite{absil08a} naturally captures such \emph{symmetries} in the search space and \emph{submerses} them into an abstract quotient manifold, where non-uniqueness \change{of the optimal solution disappears}. In essence, the framework conceptually transforms the constrained optimization problem (\ref{eq:matrix_completion}) into an \emph{unconstrained} optimization problem on the quotient manifold of fixed-rank matrices \cite{absil08a, mishra14a}. \change{An essential step of the framework is to endow} the search space with a \emph{Riemannian metric}, a smoothly varying inner product. The metric structure plays a pivotal role in relating abstract notions on quotient manifolds to concrete matrix representations. Additionally, the performance of a gradient scheme on a manifold is \emph{profoundly} dependent on the chosen metric \cite[Chapter~3]{nocedal06a}. A potential limitation of the Riemannian framework is in identifying a suitable Riemannian metric \cite{manton02a}. We resolve this issue by proposing a Riemannian metric that is tailored to the cost function.



Out of many different works on low-rank matrix completion, we focus on the three recent algorithms \cite{keshavan10a, ngo12a, wen12a} that have shown better performance in a number of challenging instances. Given a rank-$r$ matrix $\mat{X} \in \mathbb{R}^{n\times m}$, \cite{keshavan10a} exploits the fixed-rank factorization $\mat{X} = \mat{US V}^T$ to embed the rank constraint, where  $\mat{U}$ and $\mat{V}$ are column-orthonormal full-rank matrices of size $n\times r$ and $m\times r$, and ${\mat S} \in \mathbb{R}^{r \times r}$. At each iteration, \cite{keshavan10a} first updates $\mat{U}$ and $\mat{V}$ on the \emph{bi-Grassmann} manifold $\Gr{r}{n} \times \Gr{r}{m}$, where $\Gr{r}{n}$ is the set of $r$-dimensional subspaces in $\mathbb{R}^n$. Subsequently, a least-squares problem is solved to update $\mat S$. Building upon \cite{keshavan10a}, Ngo and Saad \cite{ngo12a} propose a \emph{matrix scaling} on the bi-Grassmann manifold to \emph{accelerate} the algorithm of \cite{keshavan10a}. In particular, Ngo and Saad \cite{ngo12a} motivate the matrix scaling as an \emph{adaptive preconditioner} for the optimization problem (\ref{eq:matrix_completion}) and implement a conjugate-gradient algorithm. The same matrix scaling also appears in \cite{wen12a} where the authors motivate their \emph{Gauss-Seidel} algorithm on the fixed-rank matrix factorization $\mat{X} = \mat{GH}^T$, where $\mat{G}$ and $\mat H$ are full column-rank matrices of size $n \times r$ and $m \times r$, respectively. $\mat{G}$ and $\mat H$ are updated \emph{alternatively}. A potential limitation of these algorithms is that they are alternating minimization and first-order algorithms, and extending them to other classes of optimization methods is not trivial.

In order to get the best of these methods, we reinterpret the matrix scaling of \cite{wen12a, ngo12a} as the outcome of a \emph{specific} Riemannian metric construction and exploit it to propose a novel Riemannian geometry for $\mathbb{R}_r^{n\times r}$ and an algorithm for (\ref{eq:matrix_completion}). The present work \emph{follows and completes} the \emph{unpublished} technical report \cite{mishra12a}, where we propose the general idea of tuning the Riemannian metric on a particular fixed-rank factorization. The metric so constructed provides a trade-off between second-order information of the cost function, symmetries in the search space, and simplicity of the computations involved. Building upon this work, we propose a three-factor \emph{SVD-type} factorization, and exploit this for (\ref{eq:matrix_completion}).

Our contribution is the nonlinear conjugate-gradient algorithm {\rmc} (Algorithm \ref{alg:algorithm}) that is based on a novel three-factor Riemannian geometry for matrix completion proposed in Section \ref{sec:matrix_factorization}. Although the new quotient geometry enables to propose second-order methods like the Riemannian trust-region method, here we specifically focus on conjugate-gradients as they offer an appropriate balance between convergence and computational cost. They have shown superior performance in our examples. The geometric notions to implement an off-the-shelf conjugate-gradient algorithm are listed in Section \ref{sec:optimization_ingredients}. The concrete computations for the matrix completion problem are presented in Section \ref{sec:algorithm}. Finally in Section \ref{sec:numerical_comparisons}, we illustrate the performance of {\rmc} across different problem instances, focusing in particular, on scarcely sampled and ill-conditioned problems. A Matlab implementation of {\rmc} is available on \url{http://www.montefiore.ulg.ac.be/~mishra/codes/R3MC.html}. The generic implementation of the algorithm for trust-regions is supplied with the Manopt package from \url{http://manopt.org} \cite{boumal14a}.


\section{Three-factor fixed-rank matrix factorization}\label{sec:matrix_factorization}
We consider a three-factor matrix factorization of rank-$r$ matrix $\mat{X}\in \mathbb{R}_r ^{n\times m}$, $\mat{X} = \mat{URV}^T$, where $\mat{U} \in \Stiefel{r}{n}$, $\mat{V} \in \Stiefel{r}{m}$ and ${\mat R} \in \mathbb{R}_*^{r \times r}$. The Stiefel manifold $\Stiefel{r}{n}$ is the set of matrices of size $n \times r$ with orthonormal columns and $\mathbb{R}_{*}^{r\times r}$ is the set of matrices of size $r \times r $ with non-zero determinant (full-rank square matrices). This factorization owes itself to the thin singular value decomposition (SVD). The difference with respect to SVD is that for the factor $\mat R$ we relax the diagonal constraint to accommodate any full-rank square matrix. 

Consequently, the optimization problem (\ref{eq:matrix_completion}) is a problem in $(\mat{U}, \mat{R}, \mat{V}) \in \Stiefel{r}{n} \times  \mathbb{R}_*^{r \times r} \times \Stiefel{r}{m}$. However, the critical points in the space $\Stiefel{r}{n} \times  \mathbb{R}_*^{r \times r} \times \Stiefel{r}{m}$ are not isolated as we have the symmetry
\begin{equation}\label{eq:symmetry}
\mat{X} = \mat{URV}^T = (\mat{U}\mat{O}_1) \mat{O}_1^T \mat{R} \mat{O}_2  (\mat{V}\mat{O}_2),
\end{equation}
for all $\mat{O}_1, \mat{O}_2 \in \OG{r}$, the set of orthogonal matrices of size $r\times r$. In other words, the matrix $\mat{X} \in \mathbb{R}_r^{n \times m}$ (and so, the cost function) remains unchanged under the group action (\ref{eq:symmetry}). Factoring out this symmetry, the problem (\ref{eq:matrix_completion}) is \emph{strictly} an optimization problem on the set of \emph{equivalence classes} $[(\mat{U},\mat{R}, \mat{V})]  = \{ (\mat{UO}_1,\mat{O}_1^T\mat{RO}_2, \mat{VO}_2): 
 (\mat{O}_1,\mat{O}_2) \in \OG{r} \times \OG{r}\}$ rather than $\Stiefel{r}{n} \times  \mathbb{R}_*^{r \times r} \times \Stiefel{r}{m}$. Defining $\overline{\mathcal M} := \Stiefel{r}{n} \times \mathbb{R}_*^{r \times r } \times \Stiefel{r}{m}$, the set of equivalence classes is denoted as $\mathcal{M} := \overline{\mathcal M} / (\OG{r} \times \OG{r})$, where $\overline{\mathcal M}$ is called the \emph{total space} (computational space) that is the product space $\Stiefel{r}{n} \times \mathbb{R}_*^{r \times r } \times \Stiefel{r}{m}$. The set $\OG{r} \times \OG{r}$ is the called the \emph{fiber space}. The set of equivalence classes $\mathcal{M}$ has the structure of a quotient matrix manifold \cite{absil08a}. Equivalently, the set $\mathbb{R}_r^{n\times m}$ of $n\times m$ rank-$r$ matrices is identified with $\mathcal{M}$, i.e., $\mathbb{R}_r^{n\times m} \simeq \mathcal{M}$.

Our main interest in studying the three-factor factorization is that it separates the \emph{scaling} and the \emph{subspace} information of a matrix as in the SVD factorization. The scaling information of $\mat{X}$ is in $\mat R$ whereas the subspace information is contained in $(\mat{U}, \mat{V})$. This separation of scaling and subspaces leads to a robust behavior in numerical comparisons as shown later in Section \ref{sec:numerical_comparisons}. Also since the fiber space is $\OG{r} \times \OG{r}$, by fixing $\mat R$ the search space is precisely the bi-Grassmann manifold $\Gr{r}{n}\times \Gr{r}{m}$ that comes into play in \cite{keshavan10a}, where $\Gr{r}{n} := \Stiefel{r}{m}/\OG{r}$ \cite{absil08a}. In other words, the three-factor geometry proposed in this paper generalizes the bi-Grassmann geometry of \cite{ngo12a, keshavan10a}.

\textbf{A novel Riemannian metric.} We represent an element of the quotient space $\mathcal M$ by $x$ (each element is an equivalence class) and its matrix representation in the total space $\overline{\mathcal M}:=\Stiefel{r}{n} \times \mathbb{R}_*^{r \times r } \times \Stiefel{r}{m}$ is given by $\bar x = (\mat{U}, \mat{R}, \mat{V})$. Equivalently, $x = [\bar x] = [(\mat{U},\mat{R},\mat{V})]$. The tangent space $T_{\bar x} \overline{\mathcal{M}}$ of the total space at $\bar{x} \in \overline{\mathcal{M}}$ is the product space of the tangent spaces of the individual manifolds. From \cite[Example~3.5.2]{absil08a} we have the matrix characterization $ T_{\bar x} \overline{\mathcal{M}} = \{ (\mat{Z}_{\mat U}, \mat{Z}_{\mat R}, \mat{Z}_{\mat V}) \in \mathbb{R}^{n\times r} \times \mathbb{R}^{r\times r} \times \mathbb{R}^{m\times r}: \mat{U}^T \mat{Z}_{\mat U} + \mat{Z}_{\mat U}^T \mat{U} = \mat{0} , 
 \mat{V}^T \mat{Z}_{\mat V} + \mat{Z}_{\mat V}^T \mat{V} = \mat{0} 
\}$.

The abstract quotient search space $\mathcal M$ is given the structure of a Riemannian quotient manifold by choosing a Riemannian metric in $\overline{\mathcal M}$ \cite{absil08a}. The metric defines an inner product between tangent vectors on the tangent space of $\overline{\mathcal{M}}$. Because the total space is a product space of well-studied manifolds, a typical metric on the total space is derived from choosing the natural metric of the product space, e.g., combining the natural metrics for $\Stiefel{r}{n}$ and $\mathbb{R}_*^{r\times r}$ \cite[Example~3.6.4]{absil08a}. However, this is not the only choice. Here we derive a different metric from the Hessian information of the cost function. The full Hessian information of (\ref{eq:matrix_completion}) is computationally costly \cite{buchanan05a}. To circumvent the issue, we consider a simplified (but related) version of the cost function in (\ref{eq:matrix_completion}). Specifically, consider the least-squares cost function $\| \mat{X} - \mat{X}^\star  \|_F ^2$ that is a simplification of the cost function in (\ref{eq:matrix_completion}) by assuming that $\Omega$ contains the full set of indices. It should be stressed that the cost function $\| \mat{X} - \mat{X}^\star  \|_F ^2$ is \emph{convex and quadratic} in $\mat{X}$, therefore, also convex and quadratic in each of the arguments $(\mat{U}, \mat{R}, \mat{V})$ \emph{individually}. Additionally, the orthogonality constraint $\Stiefel{r}{n}$ is also quadratic. This particular structure, \emph{individually-convex quadratic optimization on quadratic constraints}, plays a critical role. For example, the second-order information of the cost function becomes a relevant source of second-order information of the optimization problem \cite[Chapter~18]{nocedal06a}. Specifically, the \emph{block approximation} of the Hessian of $\| \mat{URV}^T - \mat{X}^\star  \|_F ^2$ has the matrix characterization $((\mat{RR}^T)\otimes \mat{I}_n, \mat{I}_r, (\mat{R}^T\mat{R})\otimes \mat{I}_m)$, where $\mat{I}_n$ is an $n \times n$ identity matrix, $\otimes$ is the Kronecker product of matrices, and $(\mat{U}, \mat{R}, \mat{V}) \in \overline{\mathcal M}$. The block approximation, equivalently, induces the metric
\begin{equation}\label{eq:metric}
\begin{array}{lll}
\bar{g}_{\bar x}(\bar{\xi}_{\bar x}, {\bar \eta}_{\bar x}) &=& \trace((\mat{RR}^T) {\bar \xi}_{\mat U}^T {\bar \eta}_{\mat U}) + \trace({\bar \xi}_{\mat R}^T {\bar \eta}_{\mat R})\\
&  +& \trace((\mat{R}^T \mat{R}) {\bar \xi}_{\mat V}^T {\bar \eta}_{\mat V})
\end{array}
\end{equation}
on $T_{\bar x}\overline{\mathcal{M}}$, where $\bar{\xi}_{\bar x}, \bar{\eta}_{\bar x} \in T_{\bar x}  \overline{\mathcal{M}}$ are tangent vectors with matrix characterizations $(\bar{\xi}_{\mat U}, \bar{\xi}_{\mat R}, \bar{\xi}_{\mat V})$ and $(\bar{\eta}_{\mat U}, \bar{\eta}_{\mat R}, \bar{\eta}_{\mat V})$, respectively. Since the function $\| \mat{URV}^T - \mat{X}^\star  \|_F ^2$ is constant under the action (\ref{eq:symmetry}), the $\bar{g}_{\bar x}$ induced from the block approximation of the Hessian of $\| \mat{URV}^T - \mat{X}^\star  \|_F ^2$ respects the symmetry and is a valid Riemannian metric candidate on the total space $\overline{\mathcal M}$ \cite[Section~3.6.2]{absil08a}. This is also equivalent to the metric proposed in \cite{mishra12a} (for a different fixed-rank factorization).


\section{Optimization-related ingredients}\label{sec:optimization_ingredients}
Once the metric is fixed on $\overline{\mathcal M} = \Stiefel{r}{n} \times \mathbb{R}_*^{r \times r } \times \Stiefel{r}{m}$, the concrete matrix representations necessary for implementing an off-the-shelf nonlinear conjugate-gradient algorithm for any smooth cost function on $\mathcal{M} = \overline{\mathcal M}/ (\OG{r} \times \OG{r})$ are listed in this section. These are derived systematically using principles laid down in \cite{absil08a}, the basis of which is the theory of \emph{Riemannian submersion}.

\textbf{Horizontal space.} The matrix representation of the tangent space $T_x \mathcal{M}$ of the abstract quotient manifold $\mathcal{M} = \overline{\mathcal M} / (\OG{r} \times \OG{r})$ is identified with a subspace of the tangent space of the total space $T_{\bar x} \overline{\mathcal M}$ that does not produce a displacement along the equivalence classes. This subspace is called the \emph{horizontal space} $\mathcal{H}_{\bar x} $ \cite[Section~3.6.2]{absil08a}. In particular, we decompose $T_{\bar x} \overline{\mathcal M}$ into two \emph{complementary subspaces} as $T_{\bar x} \overline{\mathcal M} = \mathcal{H}_{\bar x}  \oplus \mathcal{V}_{\bar x} $. The \emph{vertical space} $\mathcal{V}_{\bar x}$ is the tangent space to the equivalence class $[\bar x]$ and has the matrix expression $(\mat{U\Omega}_1,\mat{R\Omega}_2 - {\mat \Omega}_1 \mat{R}, \mat{V\Omega}_2 )$, where ${\mat \Omega}_1$ and ${\mat \Omega}_2$ are any skew-symmetric matrices of size $r\times r$ \cite[Example~3.5.3]{absil08a}. Its complementary subspace with respect to the metric ${\bar g}_{\bar x}$ is chosen as the horizontal space $\mathcal{H}_{\bar x} $. The horizontal space provides a matrix characterization to the abstract tangent space $T_x \mathcal{M}$. After a routine computation, the subspace has the characterization $\mathcal{H}_{\bar x}  = \{\bar{\eta}_{\bar x} \in T_{\bar x} \overline{\mathcal M}: \mat{RR}^T \mat{U}^T {\bar \eta}_{\mat U}  + {\bar \eta}_{\mat R} \mat {R}^T {\rm \ and\ } \mat{R}^T \mat{R} \mat{V}^T {\bar \eta}_{\mat V} - \mat{R}^T {\bar \eta}_{\mat R} {\rm \ are \ symmetric} \} $. Any tangent vector $\xi_{x} \in T_{x} \mathcal{M}$ (on the quotient space) is uniquely represented by a vector $\bar{\xi}_{\bar x}$ in the horizontal space $\mathcal{H}_{\bar x} $, also called its \emph{horizontal lift}. 

\textbf{Projection operators.} We require two projection operators: one from the ambient space $\mathbb{R}^{n\times r}\times \mathbb{R}^{r \times r}\times \mathbb{R}^{m \times r}$ onto the tangent space $T_{\bar x}\overline{\mathcal M}$, and a second projection operator that further extracts the horizontal component of a tangent vector. Given a matrix in the space $\mathbb{R}^{n\times r}\times \mathbb{R}^{r \times r}\times \mathbb{R}^{m \times r}$, its projection onto the tangent space $T_{\bar x} \overline{\mathcal M}$ is obtained by extracting the component normal, in the metric sense, to the tangent space. The normal space $N_{\bar x} \overline{\mathcal M}$ is, thus, $\{ (\mat{UN}_1, \emptyset, \mat{VN}_2): \mat{N}_2\mat{RR}^T{\rm \ and \ } \mat{N}_2\mat{R}^T\mat{R}{\rm \ are \  symmetric}   \}$ for $\mat{N}_1, \mat{N}_2 \in \mathbb{R}^{r \times r}$. Extracting the tangential component is accomplished by using the linear operator $\Psi_{\bar x}: \mathbb{R}^{n\times r}\times \mathbb{R}^{r \times r}\times \mathbb{R}^{m \times r} \rightarrow T_{\bar x} \overline{\mathcal M}:( \mat{Z}_{\mat U},\mat{Z}_{\mat R},  \mat{Z}_{\mat V}) \mapsto \Psi_{\bar x}( \mat{Z}_{\mat U},\mat{Z}_{\mat R},  \mat{Z}_{\mat V})$, i.e.,
\begin{equation}\label{eq:projection_tangent_space}
\begin{array}{lll}
\Psi_{\bar x}( \mat{Z}_{\mat U},\mat{Z}_{\mat R},  \mat{Z}_{\mat V}) &= & (\mat{Z}_{\mat U} - \mat{U}\mat{B}_{\mat U} (\mat{RR}^T)^{-1}, \mat{Z}_{\mat R}, \\
& & \mat{Z}_{\mat V} - \mat{V}\mat{B}_{\mat V} (\mat{R}^T\mat{R})^{-1}   ),
\end{array}
\end{equation}
where $\mat{B}_{\mat U}$ and $\mat{B}_{\mat V}$ are symmetric matrices of size $r\times r$ that are obtained by solving the Lyapunov equations
$ \mat{RR}^T \mat{B}_{\mat U} + \mat{B}_{\mat U} \mat{RR}^T = \mat{RR}^T (\mat{U}^T \mat{Z}_{\mat U} + \mat{Z}_{\mat U}^T \mat{U}        ) \mat{RR}^T $ and $\mat{R}^T\mat{R} \mat{B}_{\mat V} + \mat{B}_{\mat V} \mat{R}^T\mat{R} = \mat{R}^T\mat{R}  (\mat{V}^T \mat{Z}_{\mat V} + \mat{Z}_{\mat V}^T \mat{V}        )\mat{R}^T\mat{R} $ which can be solved efficiently with a total cost $O(r^3)$.


A subsequent projection onto the horizontal space is obtained by the operator $\Pi_{\bar x}: T_{\bar x} \overline{\mathcal M} \rightarrow \mathcal{H}_{\bar x} : \bar{\xi}_{\bar x} \mapsto \Pi_{\bar x}(\bar{\xi}_{\bar x})$
\begin{equation}\label{eq:projection_horizontal_space}
\begin{array}{lll}
\Pi_{\bar x}(\bar{\xi}_{\bar x}) =( \bar{\xi}_{\mat U} - \mat{U}\mat{\Omega}_1, \bar{\xi}_{\mat R} + \mat{\Omega}_1 \mat{R} - \mat{R}\mat{\Omega}_2 ,
     \bar{\xi}_{\mat V} - \mat{V}\mat{\Omega}_2 ),
\end{array}
\end{equation}
where $\mat{\Omega}_1$ and $\mat{\Omega}_2$ are skew-symmetric matrices of size $r\times r$ that are obtained by solving the coupled system of Lyapunov equations that has the form
\begin{equation}\label{eq:lyap_horizontal_space}
\begin{array}{lll}
\mat{RR}^T \mat{\Omega}_1 + & \mat{\Omega}_1 \mat{RR}^T - \mat{R}\mat{\Omega}_2 \mat{R}^T \\
 & = \Skew(\mat{U}^T\bar{\xi}_{\mat U} \mat{RR}^T ) + \Skew(\mat{R}\bar{\xi}_{\mat R}^T), \\
\mat{R}^T\mat{R} \mat{\Omega}_2 + & \mat{\Omega}_2 \mat{R}^T\mat{R} - \mat{R}^T\mat{\Omega}_1 \mat{R} \\
& = \Skew(\mat{V}^T\bar{\xi}_{\mat V} \mat{R}^T\mat{R} ) + \Skew(\mat{R}^T\bar{\xi}_{\mat R}), 
\end{array}
\end{equation}
where $\Skew(\cdot)$ extracts the skew-symmetric part of a square matrix, i.e., $\Skew(\mat D) = (\mat{D} - \mat{D}^T)/2$. The coupled equations (\ref{eq:lyap_horizontal_space}) can also be solved efficiently (and in closed form) by diagonalizing $\mat R$ and performing similarity transforms on ${\mat \Omega}_1$ and ${\mat \Omega}_{2}$. The computational cost is $O(r^3)$.

\textbf{Retraction.} A retraction is a mapping that maps vectors in the horizontal space to points on the search space $\overline{\mathcal{M}}$. Precise formulation is in \cite[Definition~4.1.1]{absil08a}). It provides a natural way to move on the manifold along a horizontal direction. The product nature of the total space $\overline{\mathcal M}$ again allows to choose a retraction by simply combining the retractions on the individual manifolds \cite[Example~4.1.3]{absil08a}
\begin{equation}\label{eq:retraction}
\begin{array}{ll}
R_{\bar x} (\bar{\xi}_{\bar x}) = ({\rm uf}(\mat{U} + \bar{\xi}_{\mat U}),  \mat{R} + \bar{\xi}_{\mat R} , {\rm uf}(\mat{V} + \bar{\xi}_{\mat V})),
\end{array}
\end{equation}
where $\bar{\xi}_{\bar x} \in \mathcal{H}_{\bar x} $ and ${\rm uf}(\cdot)$ extracts the orthogonal factor of a full column-rank matrix, i.e., ${\rm uf}(\mat A) = \mat{A}(\mat{A}^T \mat{A})^{-1/2}$ which is computed efficiently. The computational cost of a retraction operation is $O(nr^2 + mr^2 + r^3)$.

\textbf{Vector transport.} A vector transport $\mathcal{T}_{\eta_x} \xi_x$ on a manifold $\mathcal{M}$ is a smooth mapping that transports a tangent vector $\xi_x \in T_x \mathcal{M}$ at $x \in \mathcal{M}$ to a vector in the tangent space at $R_{x}(\eta _x)$ satisfying the first-order approximation of transporting a vector along the geodesic \cite[Definition~8.1.1]{absil08a}. With the projection operators, defined earlier, the horizontal lift of the vector transport $\mathcal{T}_{\eta_x} \xi_x$ is 
\begin{equation}\label{eq:vector_transport}
\overline{\mathcal{T}_{\eta_x} \xi_x} =  \Pi_{R_{\bar x}(\bar{\eta} _{\bar x})} (\Psi_{R_{\bar x}(\bar{\eta} _{\bar x})} ({\bar \xi}_{\bar x})),
\end{equation}
where $\bar{\eta}_{\bar x}$ and $\bar{\xi}_{\bar x}$ are the horizontal lifts of $\xi _x$ and $\eta _x$; and $\Pi(\cdot)$ and $\Psi(\cdot)$ are the projection operators defined in (\ref{eq:projection_horizontal_space}) and (\ref{eq:projection_tangent_space}), respectively. The computational cost of transporting a vector solely depends on projection operations, defined earlier, which cost $O(nr^2 + mr^2 + r^3)$.



\section{{\rmc}: Algorithmic details}\label{sec:algorithm}
In Section \ref{sec:optimization_ingredients}, the matrix representations of various notions on $\mathcal{M} = \overline{\mathcal M}/(\OG{r} \times \OG{r})$ are presented. Here we specifically deal with (\ref{eq:matrix_completion}) that is, equivalently, the problem
\[
\begin{array}{lll}
\min\limits_{x \in \mathcal{M}} f(x), 
\end{array}
\]
where $f: \mathcal{M} \rightarrow \mathbb{R}$ is the restriction of the function $\bar{f}: \overline{\mathcal{M}} \rightarrow \mathbb{R} : \bar{x} \mapsto \bar{f}(\bar x) = \|\mathcal{P}_{\Omega}(\mat{URV}^T) - \mathcal{P}_{\Omega}(\mat{X}^{\star})\|_F^2/|\Omega|$, $\bar{x} \in \overline{\mathcal M}$ has the matrix representation $(\mat{U}, \mat{R}, \mat{V}) \in \Stiefel{r}{n} \times \mathbb{R}_*^{r\times r} \times
\Stiefel{r}{m}$. The steps of our nonlinear conjugate-gradient algorithm {\rmc} are shown in Algorithm \ref{alg:algorithm}. The computational cost per iteration of {\rmc} is $O(|\Omega|r + nr^2 + mr^2 + r^3)$, where $|\Omega|$ is the number of known entries. The convergence of a Riemannian conjugate-gradient algorithm follows from the analysis in \cite{ring12a}. Step \ref{step:gradient_related} of Algorithm \ref{alg:algorithm}, in particular, ensures that the sequence $\{ \bar{\eta}_{i}\}$, $\bar{\eta}_i \in \mathcal{H}_{\bar{x}_i}$ is \emph{gradient-related} \cite[Definition~4.2.1]{absil08a}.

\begin{algorithm}
 \caption{{\rmc}: a conjugate-gradient algorithm for (\ref{eq:matrix_completion})}
 \label{alg:algorithm}
\begin{algorithmic}[1]
\Input Given $\bar{x}_0 = (\mat{U}_0, \mat{R}_0, \mat{V}_0)\in \overline{\mathcal M}$ and $\bar{\eta}_0 = 0$.
%
 \State Compute the gradient $\bar{\xi}_i = \overline{\grad}_{{\bar x}_{i}} \bar{f}$. 
 \Comment{(\ref{eq:Riemannian_gradient})}
\State Compute the conjugate  direction by Polak-Ribi\`ere (PR+)
\Statex $\bar{\eta}_{i}  = -{{\bar \xi}_i} + \beta _i \Pi_{\bar{x}_i}(\Psi_{\bar{x}_i} (\bar{\eta} _{i -1 }) )$. \Comment{(\ref{eq:vector_transport}) and \cite[(8.28)]{absil08a}}
\State \label{step:gradient_related} Check whether $\bar{\eta}_i$ is a descent direction. If not, then
\Statex $\bar{\eta}_i = - \overline{\grad}_{{\bar x}_i} \bar{f}$.
\State Determine an initial step $s_i$.
\Comment{Section \ref{sec:algorithm}}
\State Retract along $\bar{\eta} _i$ to compute the new iterate $\bar{x}_{i+1} $.
\Comment{(\ref{eq:retraction})}
\State Repeat until convergence.
 \end{algorithmic}
\end{algorithm}



\textbf{The gradient computation.} We define an auxiliary sparse variable $\mat{S} = 2(\mathcal{P}_{\Omega}(\mat{URV}^T) - \mathcal{P}_{\Omega}(\mat{X}^{\star}))/|\Omega| $ that is interpreted as the Euclidean gradient of $\bar{f}$ in the space $\mathbb{R}^{n \times m}$. The horizontal lift (matrix representation) of the Riemannian gradient $\grad_{x} f$ is characterized using the projection operator (\ref{eq:projection_tangent_space}) \cite[Section~3.6]{absil08a} as
\begin{equation}\label{eq:Riemannian_gradient}
\begin{array}{lll}
 \overline{\grad_{x}f} =(\mat{SVR}^T (\mat{RR}^T)^{-1} - \mat{U}\mat{B}_{\mat U} (\mat{RR}^T)^{-1},
  \mat{U}^T \mat{SV}, \\
\quad \qquad \qquad \mat{S}^T \mat{UR} (\mat{R}^T \mat{R})^{-1} - \mat{V}\mat{B}_{\mat V}(\mat{R}^T \mat{R})^{-1} ), \\
\end{array}
\end{equation}
where $\mat{B}_{\mat U}$ and $\mat{B}_{\mat V}$ are the solutions to the Lyapunov equations $\mat{RR}^T \mat{B}_{\mat U} + \mat{B}_{\mat U} \mat{RR}^T = 2\Sym(\mat{RR}^T \mat{U}^T \mat{SVR}^T)$ and 
$\mat{R}^T\mat{R} \mat{B}_{\mat V} + \mat{B}_{\mat V} \mat{R}^T\mat{R} =2\Sym( \mat{R}^T\mat{R}  \mat{V}^T \mat{S}^T\mat{UR} )$.
Here $\Sym(\cdot)$ extracts the symmetric part of a square matrix, i.e., $\Sym(\mat{D}) = (\mat{D} + \mat{D}^T)/2$. The total numerical cost of computing the Riemannian gradient is $O(|\Omega |r)$.

\textbf{Initial guess for the step-size.} The least-squares nature of the cost function in (\ref{eq:matrix_completion}) also allows to compute a \emph{linearized} step-size guess efficiently \cite{vandereycken13a, mishra12a}. Given a search direction $\bar{\xi}_{\bar x} \in \mathcal{H}_{\bar x}$, we solve the following optimization problem to obtain a linearized step-size guess
by considering a degree $2$ polynomial approximation, i.e., 
$ s_* =  \argmin _{s \in \mathbb{R}_+} \| \mathcal{P}_{\Omega}(\mat{URV}^T + s\bar{\xi}_{\mat U}\mat{RV}^T  + s\mat{U} \bar{\xi}_{\mat R} \mat{V}^T     + s\mat{UR}\bar{\xi}_{\mat V}^T )  - \mathcal{P}_{\Omega}(\mat{X}^{\star})  \|_F^2 
$ that has a closed form solution. Computing $s_*$ incurs a numerical cost $O(|\Omega|r)$.

\textbf{Rank updating.}  In many problems, a good rank of the solution is either not known a priori or the notion of numerical rank is too vague to define it precisely, e.g., matrices with exponential decay of singular values. In such instances, it makes sense to traverse through a number of ranks, and not just one, in a systematic manner while ensuring that the cost function of (\ref{eq:matrix_completion}) decreases with each such rank update. One way is to use fixed-rank optimization in conjunction with a rank-update strategy \cite{burer03a, mishra13a}. The rank-one update is based on the idea of moving along the dominant rank-one projection of the negative gradient in the space $\mathbb{R}^{n \times m}$. If $(\mat{U}, \mat{R}, \mat{V})$ is a rank-$r$ matrix, then the rank-one update corresponds to $(\mat{U}_+ ,\mat{R}_+,\mat{V}_+)$ such that $\mat{U}_+ \mat{R}_+\mat{V}_+^T = \mat{U} \mat{R} \mat{V} ^T  - \sigma u v^T$ where $u\in \mathbb{R}^n$ and $v \in \mathbb{R}^m$ are the unit-norm dominant left and right singular vectors of $\mat S$ and $\sigma > 0$ is the dominant singular value. The total computational cost is $O(|\Omega|  + nr^2 + m r^2 + r^3)$.

\section{Numerical comparisons}\label{sec:numerical_comparisons}
We compare {\rmc} with the state-of-the-art algorithms RTRMC \cite{boumal11a}, LMaFit \cite{wen12a}, {\rtwomc} (also known as qGeomMC)  \cite{mishra12a}, ScGrassMC \cite{ngo12a}, LRGeom \cite{vandereycken13a}, and Polar Factorization \cite{meyer11a, mishra14a}. For the last four algorithms, we use with their conjugate-gradient implementations. {\rmc} and Polar Factorization exploit accelerated step-size computation in Section \ref{sec:algorithm}. The choice of these algorithms as state-of-the-art rests on recent publications \cite{meyer11a, boumal11a, ngo12a, wen12a,  mishra14a, vandereycken13a} and references therein. Other problems formulations, e.g., with nuclear norm regularization \cite{cai10a}, have not been discussed here as the main motivation of the present paper is to look specifically look at algorithms that exploit the search space of fixed-rank matrices. In this section, we first show the connection of {\rmc} to ScGrassMC \cite{ngo12a} and then provide a comparison of the considered algorithms across different problem instances, including one on a real dataset.

All simulations are performed in Matlab on a $2.53$ GHz Intel Core i$5$ machine with $4$ GB of RAM. We use Matlab codes for all the  considered algorithms. For each example, an $n \times m$ random matrix of rank $r$ is generated as in \cite{cai10a}. Two matrices $\mat{A} \in \mathbb{R}^{n \times r}$ and $\mat{B} \in \mathbb{R}^{m \times r}$ are generated according to a Gaussian distribution with zero mean and unit standard deviation. The matrix product $\mat{AB} ^T$ gives a random rank-$r$ matrix. A fraction of the entries are randomly removed with uniform probability. We use the over-sampling ratio (OS), ${\rm OS} = |\Omega|/(nr +mr -r^2)$, to quantify the fraction of known entries. The maximum number of iterations of all is set to $500$ (for RTRMC, it is $200$ ). The algorithms are initialized similarly and stopped when the cost function is below $10^{-20}$.


\textbf{Connection to ScGrassMC.} The connection with ScGrassMC is apparent from the fact that it scales the Riemannian gradient on $\Gr{r}{n} \times \Gr{r}{m}$ by $((\mat{RR}^T)^{-1}, (\mat{R}^T \mat{R})^{-1} )$. It should be noted that this is the same scaling that is obtained in the gradient computation (\ref{eq:Riemannian_gradient}). Additionally, we have linear projections to respect the quotient nature of the search space of fixed-rank matrices. The difference with respect to ScGrassMC is on two fronts. First, we perform a simultaneous update of the variables $(\mat{U}, \mat{R}, \mat{V})$, while ScGrassMC alternates between updating $(\mat{U}, \mat{V})$ and $\mat R$. Second, while preconditioning is motivated in \cite{ngo12a} as a way to accelerate the algorithm of \cite{keshavan10a}, we view it as a particular choice of the Riemannian metric, enabling to develop arbitrary unconstrained optimization algorithms. 

\begin{figure*}[t]
\center
\hspace*{-1.5em}
\hspace*{-.5em}
\subfigure[${\rm OS} = 2.1$]{\label{fig:low_sampling}
\includegraphics[scale = 0.17]{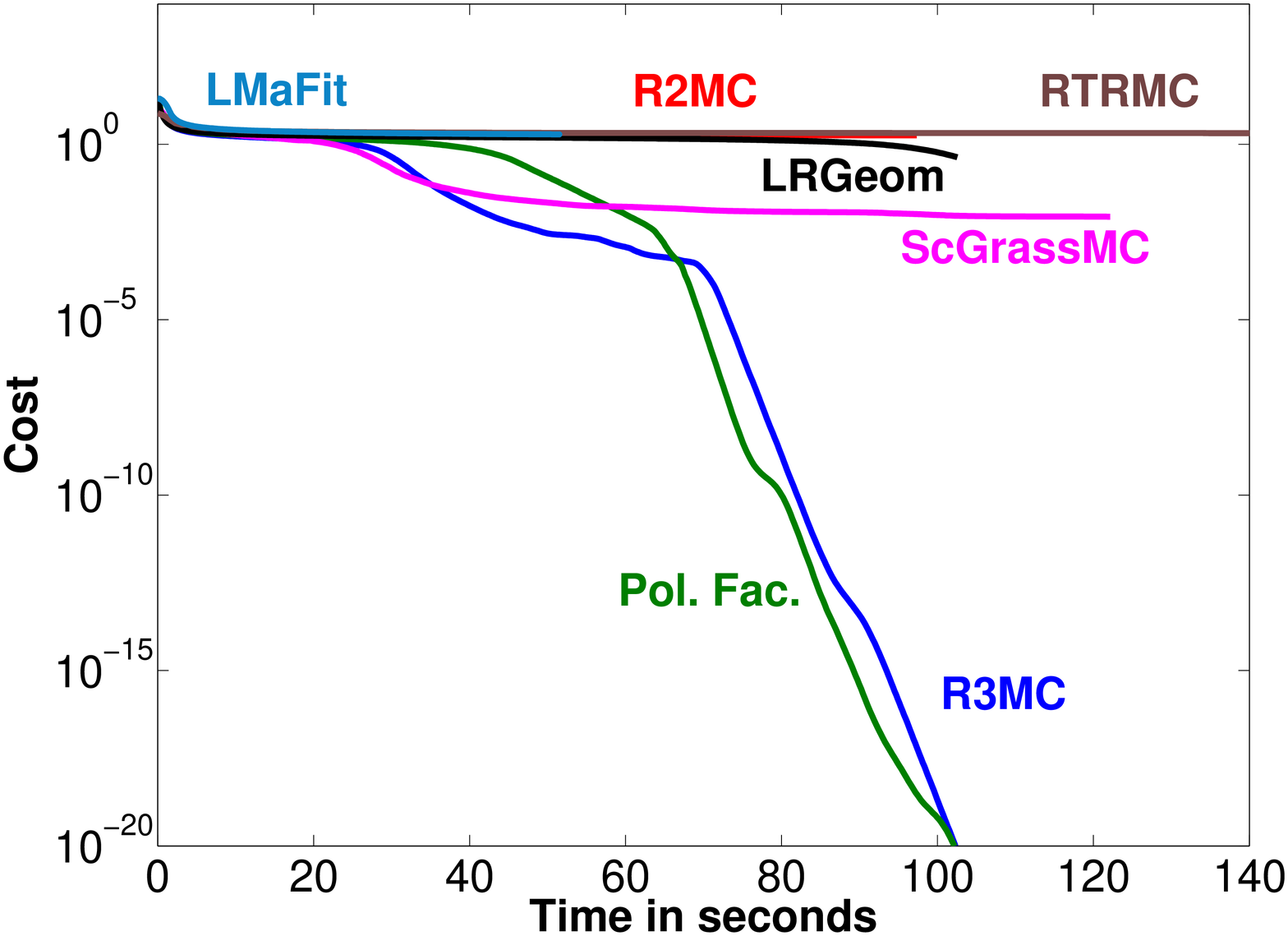}
}
\hspace*{-1.5em}
\subfigure[${\rm CN} = 500$,  ${\rm OS}=3$]{\label{fig:ill_conditioning}
\includegraphics[scale = 0.17]{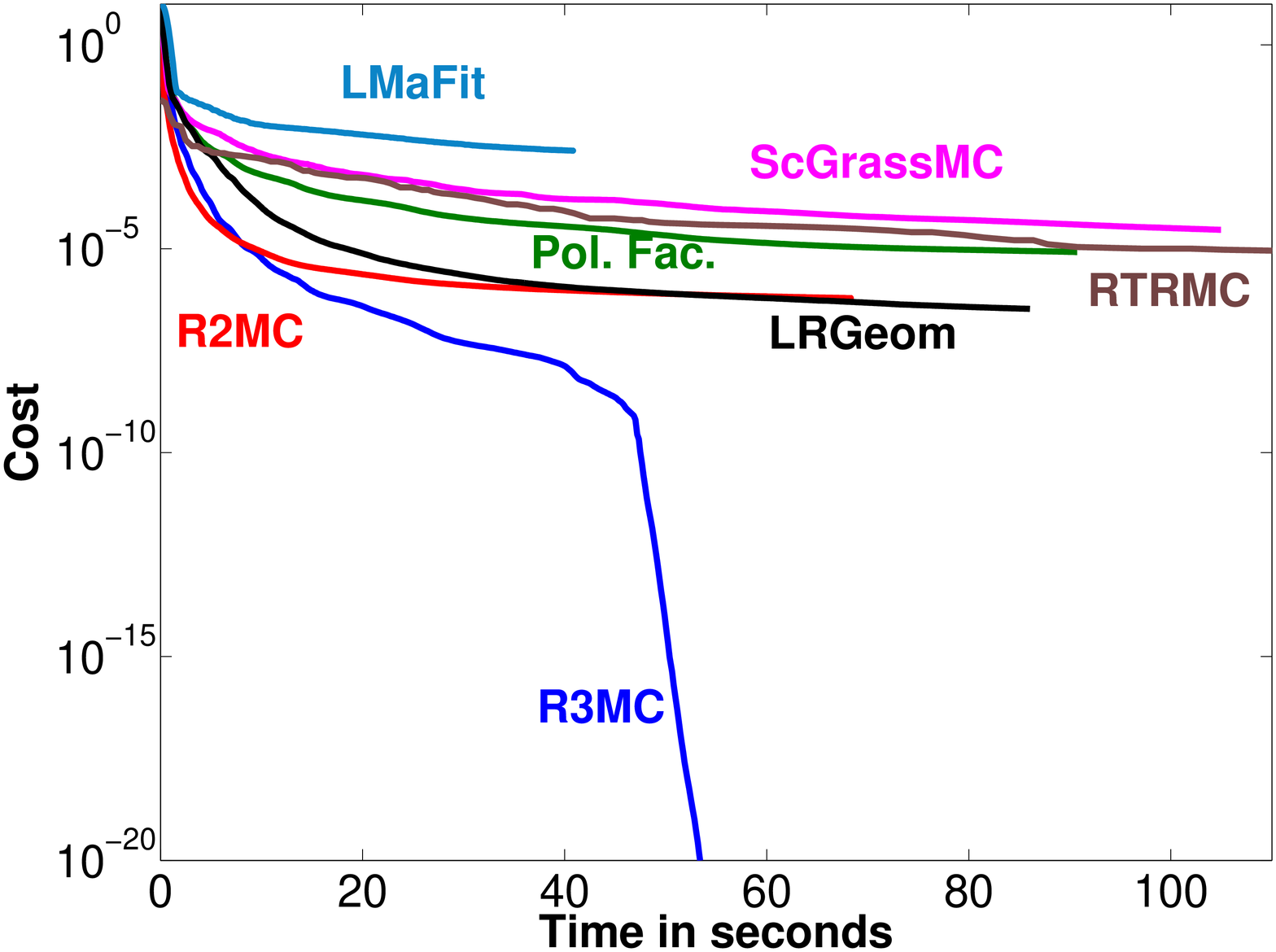}
}
\hspace*{-1.5em}
\subfigure[$\rm{OS} = 4$, ${\rm CN} = 500$]{\label{fig:low_sampling_ill_conditioning}
\includegraphics[scale = 0.17]{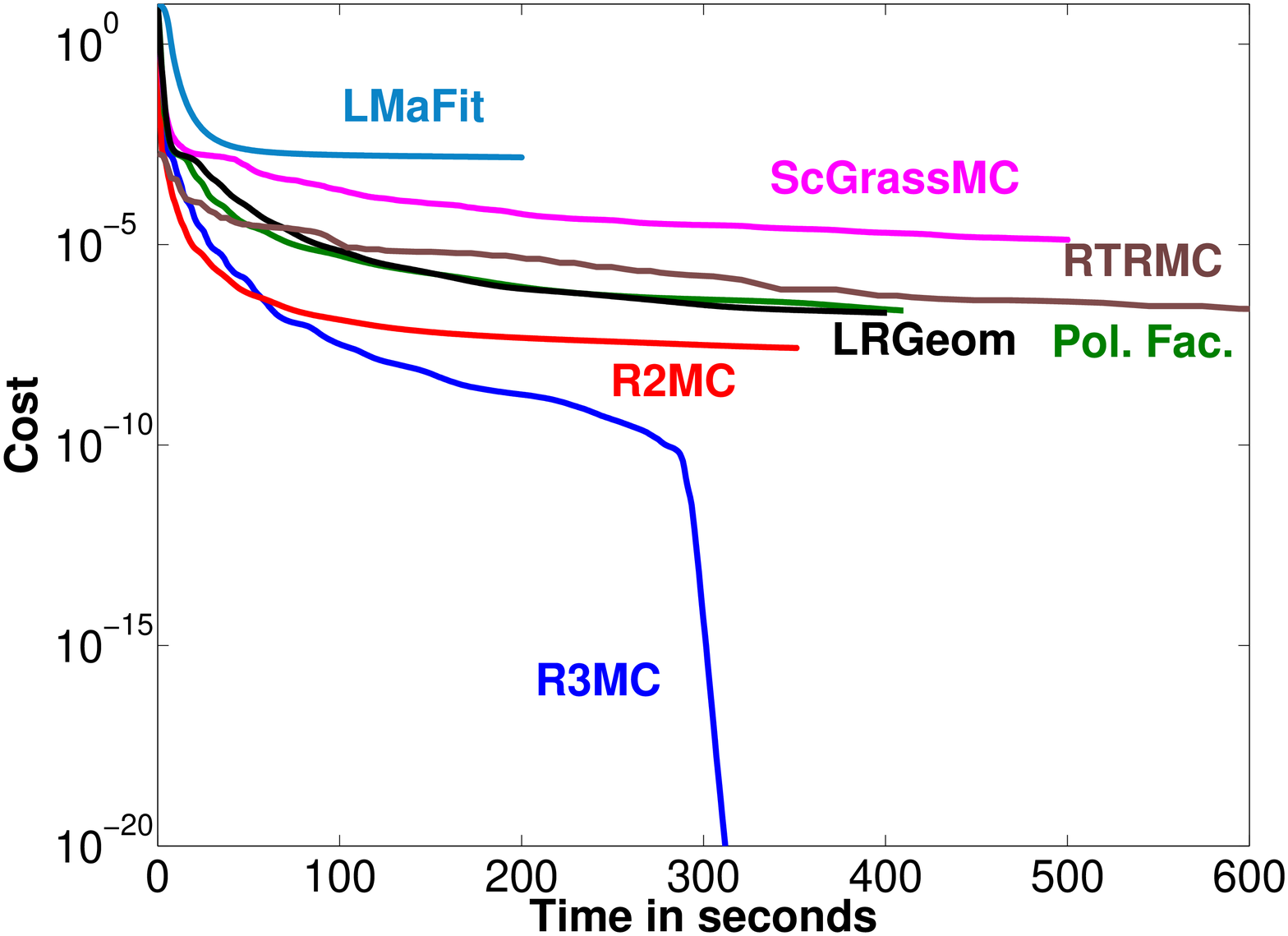}
}
\hspace*{-1.5em}
\subfigure[Rank-one updating]
{\label{fig:homotopy}
\includegraphics[scale = 0.17]{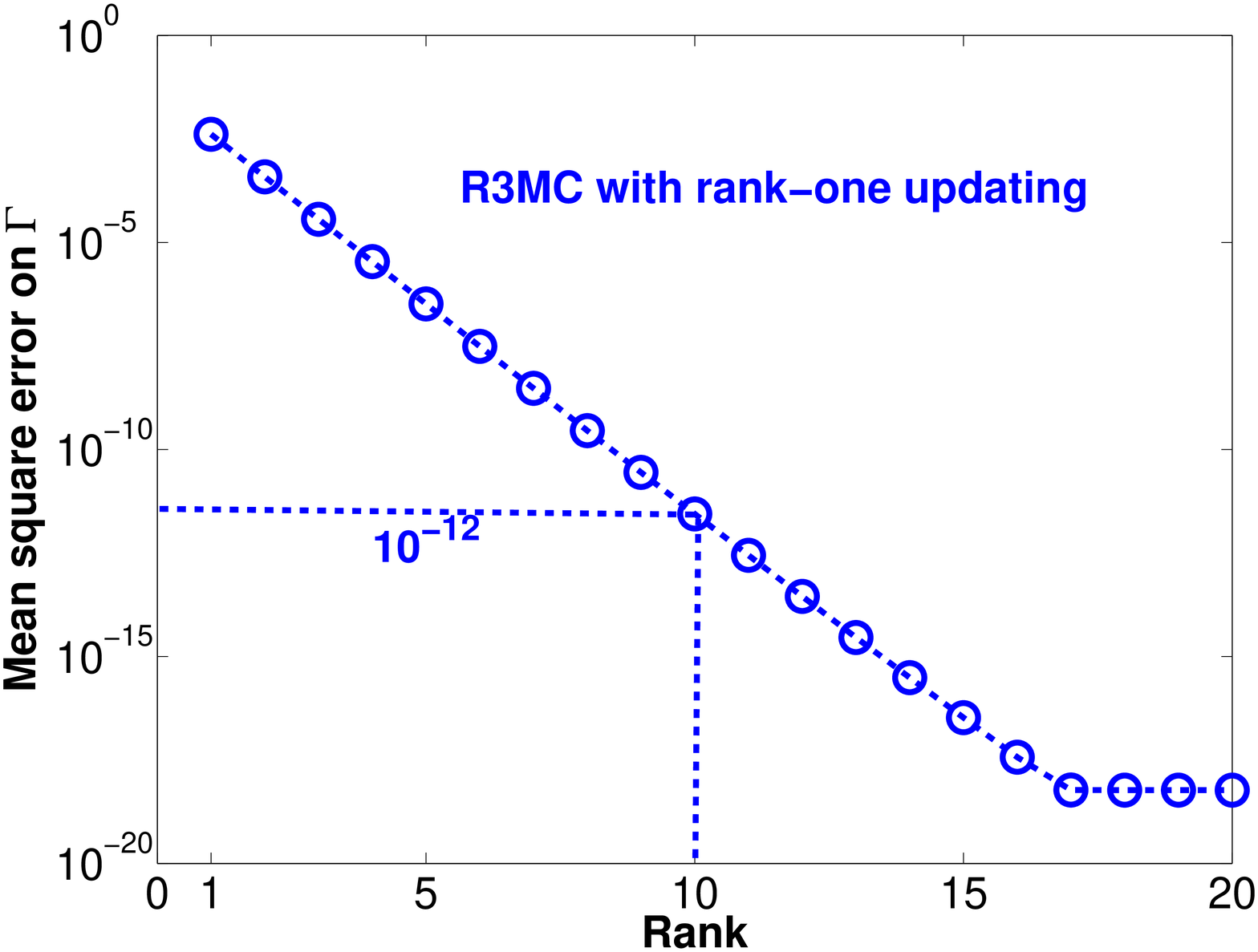}
}
\hspace*{-1.5em}
\subfigure[Rectangular matrices]{\label{fig:rectangular}
\label{fig:rectangular:test_error}
\includegraphics[scale = 0.17]{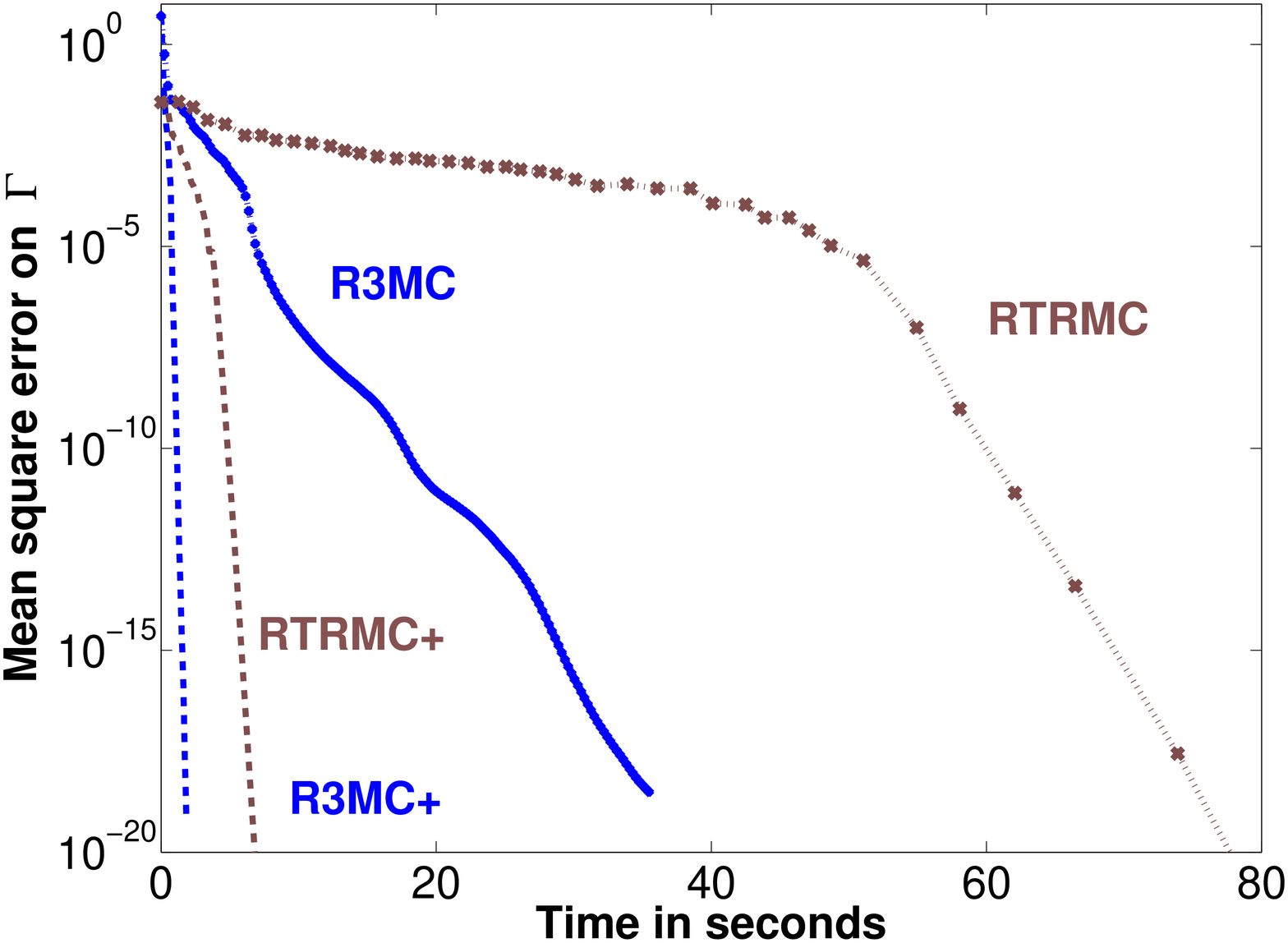}
}
\hspace*{-1.5em}
\subfigure[MovieLens dataset]{
\label{fig:movielens}
\includegraphics[scale = 0.17]{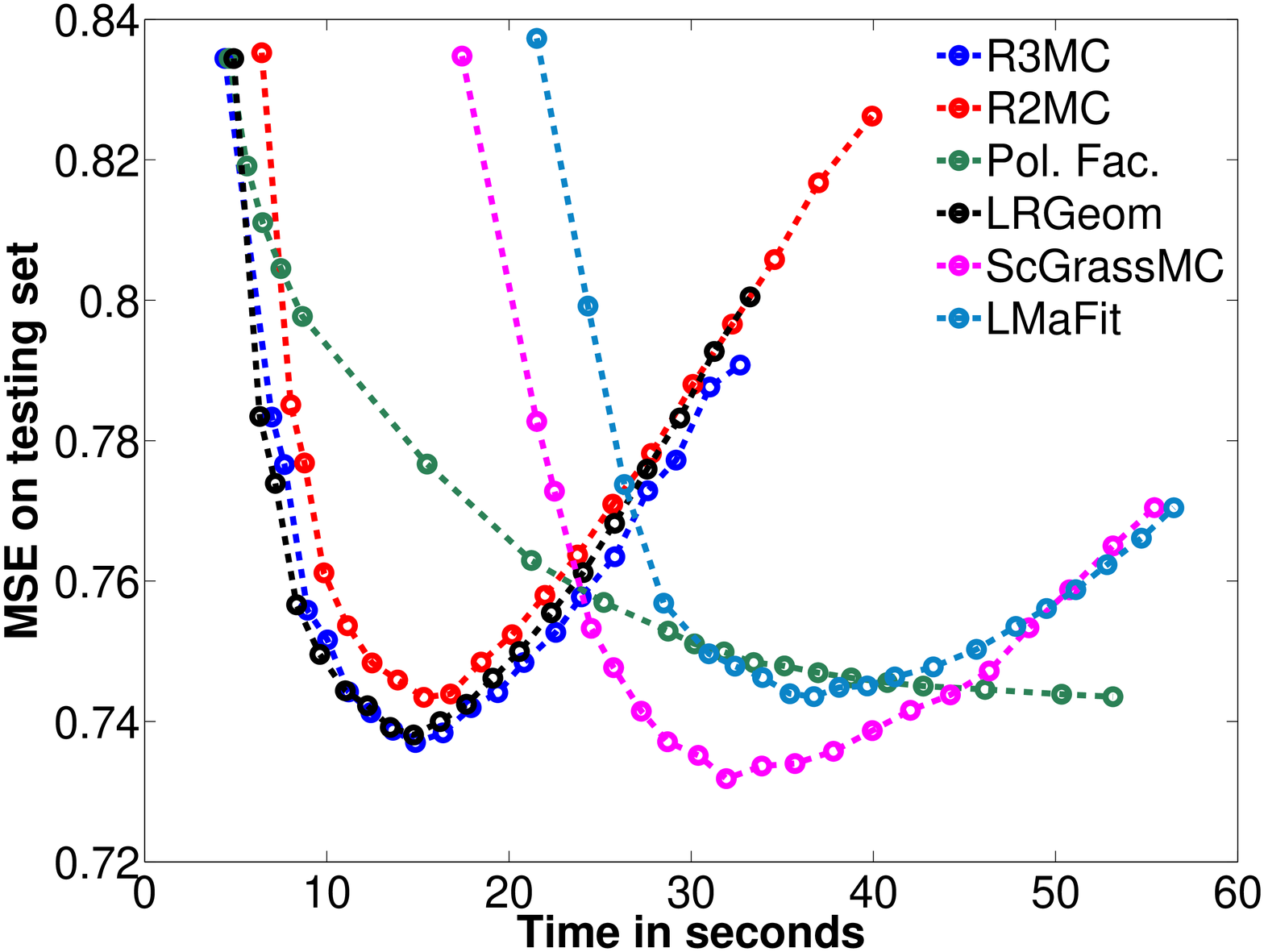}
}
\caption{Performance of different algorithms under different scenarios. {\rmc} is particularly efficient in a number of instances.}
\label{fig:overall}
\end{figure*}
\begin{table*}[t]
\caption{Mean square errors obtained on the test set of the MovieLens-$1$M dataset.} 
\begin{center} \small 
\begin{tabular}{ | p{1.4cm} | p{1.2cm}| p{1.2cm} | p{1.6cm} | p{1.4cm}| p{1.6cm}| p{1.4cm}| p{1.2cm}| } 
\hline
 &  &  &  &  &  &  & \\

$r$  & {\rmc}& {\rtwomc} & Pol. Fac. & ScGrassMC & LRGeom & LMaFit & RTRMC  \\ 
 &  &  &  &  &  &  & \\
\hline
 &  &  &  &  &  &  & \\
$3$&  $0.7713$   & $0.7771$  &  $0.7710$  & $0.7967$  & $ 0.7723$ & $0.7762$ & $0.7858$ \\ 
$4$&  $0.7677$   & $0.7758$  &  $0.7675$  & $0.7730$  & $ 0.7689$ & $0.7727$ & $0.8022$ \\ 
$5$ &$0.7666 $ & $ 0.7781$ & $0.7850 $  & $0.8280 $  & $ 0.7660$  & $ 0.8224$  & $0.8314$ \\ 
$ 6$ & $ \textbf{0.7634}$  & $0.7893$  & $0.7651$  & $ 0.7910$  & $ 0.7698$ & $ 0.8194$  & $ 0.8802$ \\ 
$7$ & $0.7684$  &   $  0.7996 $ &   $  0.7980$ &  $  0.8368 $ &  $  0.7810 $    &  $0.8074$ (max. iters)   &     $ 0.8241$ \\ 
 &  &  &  &  &  &  & \\

\hline
 &  &  &  &  &  &  & \\

With rank updates& $0.7370$ (9) & $0.7434 $ (8) & $0.7435 $ (20 max. rank) & $\textbf{0.7323} $ (10) & $ 0.7381$ (9)  & $0.7435$ (9) & - \\ 
 &  &  &  &  &  &  & \\

\hline
\end{tabular}
\end{center} 
\label{tab:movielens} 
\end{table*}

\textbf{Case (a): influence of over-sampling.} We consider a moderate scale matrix of size $10000\times 10000$ of rank $10$. Different instances with different over-sampling have been considered. For larger values of OS, most of the algorithms perform similarly and show a nice behavior. With smaller OS ratios, the algorithms, however, perform very differently. In fact in Figure \ref{fig:low_sampling}, for the case of ${\rm OS} = 2.1 $ only {\rmc} and Polar Factorization \cite{mishra14a} algorithms converged, pointing towards a robustness of the three-factor matrix factorizations.

\textbf{Case (b): influence of conditioning.} We consider matrices of size $5000\times 5000$ of rank $10$, ${\rm CN} = 3$, and impose an exponential decay of singular values. The ratio of the largest to the lowest singular value is known as the condition number (CN) of the matrix. At rank $10$ the singular values with condition number $100$ is obtained using the Matlab function \verb+logspace(-2,0,10)+. The over-sampling ratio for these instances is $3$. The matrix completion problem becomes challenging as the ${\rm CN}$ increases. Figure \ref{fig:ill_conditioning} shows a particular instance. In general, for smaller ${\rm CN}$, {\rtwomc}, {\rmc}, and LRGeom perform better than the others. In instances of larger $\rm{CN}$, {\rmc} outperforms the others.

\textbf{Case (c): influence of low sampling and ill-conditioning.} In this test, we look into problem instances which result from both scarcely sampled and ill-conditioned data. The test requires completing relatively large matrices of size $25000\times 25000$ of rank $10$ with different condition numbers and ${\rm OS}$ ratios. Figure \ref{fig:low_sampling_ill_conditioning} shows the good performance of {\rmc}.

\textbf{Case (d): ill-conditioning and rank-one updates.} We generate a random matrix of size $5000 \times 5000$ of rank $20$ with exponentially decaying singular values so that the condition number is $10^{10}$ and ${\rm OS}$ is $2$ (computed for rank $10$). Figure \ref{fig:homotopy} shows the recovery of a set of entries in $\Gamma$, that is different from $\Omega$ on which we optimize, when {\rmc} is used with the rank-one updating procedure of Section \ref{sec:algorithm}. Most fixed-rank algorithms show a better performance when combined with a rank incrementing strategy.


\textbf{Case (e): rectangular matrices.} Here we are particularly interested in instances with $n \ll m$, i.e., \emph{rectangular} matrices. For these instances, most simulations suggest that the algorithm RTRMC of \cite{boumal11a} performs numerically very efficiently. This is not surprising as the underlying geometry of RTRMC exploits the fact that the least-squares formulation (\ref{eq:matrix_completion}) is solvable by fixing one of the fixed-rank factors. 

To adapt the algorithms like {\rmc} (including RTRMC) to rectangular matrices under the standard assumptions for the matrix completion problem, we propose to deal with smaller size matrices with fewer columns. Consider a \emph{truncated} submatrix of size $n\times p$ that picks all the rows of $\mat{X}^\star$ but picks only $p$ (randomly chosen out of $m$) columns. A simple analysis shows that ${\rm OS}_{\rm trunc.} = {\rm OS}{\alpha}/({1+\alpha})$, where ${\rm OS}$ and ${\rm OS}_{\rm trunc.}$ are the over-sampling ratios for full and truncated problems, and $\alpha =  p/n $. This means that the truncated problem is \emph{challenging} for smaller values $\alpha$. However for $\alpha > 1$, it is possible to have a competitive \emph{trade-off} between difficulty and computational cost. It should, however, be noted that for a sufficiently large $\alpha$ both the original and the truncated matrices share the same \emph{left subspace} \cite{balzano10a}. Accordingly, we deal with the truncated problem of dimension $n \times p$ to compute the left subspace of the original matrix. Once the left subspace $\mat{U} \in \Stiefel{r}{n}$ is identified, the weighting factor, e.g., the matrix $\mat{W} \in \mathbb{R}^{r \times n}$ of the factorization $\mat{X} = \mat{UW}$ of \cite{boumal11a} is obtained by solving a least-squares problem by fixing $\mat{U}$ \cite{boumal11a}. A QR factorization of $\mat{W}$ results in the factors $\mat{R}$ and $\mat{V}$. The resulting $(\mat{U}, \mat{R}, \mat{V})$ provides a good initialization to algorithms for the original problem. 

We consider a rank-$5$ matrix of size $1000 \times 50000$ with ${\rm OS} = 5$ and ${\rm CN} = 10$. An incomplete submatrix of size $1000 \times 2000$ is formed by picking randomly $2000$ columns. Consequently, $\alpha = 2$ and $ {\rm OS}_{\rm trunc.} = 3.3$. The mean square errors on a set of entries $\Gamma$, different from $\Omega$, are reported in Figure \ref{fig:rectangular}, where both {\rmc} and RTRMC with the proposed scheme (appended by the sign +) are \emph{significantly} faster than their counterparts that deal with the full problem.

\textbf{Case (f): MovieLens dataset.} As a final test, we compare the algorithms on the MovieLens-$1$M dataset, downloaded from \url{http://grouplens.org/datasets/movielens/}. The dataset has a million ratings corresponding to $6040$ users and $3952$ movies. We perform $10$ random $80/10/10$-train/validation/test partitions of the ratings. The algorithms are run on the train set. The results are reported on the test set, averaged over $10$ partitions. In RTRMC, the parameter $\lambda$ is set to $ 10^{-6}$, to avoid the error due to non-uniqueness of the least-squares solution that it uses. Finally, the maximum number of iterations is set to $1000$ ($200$ for RTRMC). The algorithms are stopped before when the mean square error (MSE) on the validation set starts to increase. Table \ref{tab:movielens} (the second row) shows the MSEs on the test set with standard deviations $\pm 10^{-5}$ for different ranks. The best score of $0.7634$ is obtained by {\rmc} at rank $6$. Figure \ref{fig:movielens} shows the time taken by the algorithms, where {\rmc}, {\rtwomc}, and LRGeom are faster than the others. We also run all the algorithms with the rank-updating procedure of Section \ref{sec:algorithm} by traversing through all the ranks from $1$ to $20$. The rank is updated when the error on the validation set starts to increase. The last row of Table \ref{tab:movielens} compiles the \emph{best} MSEs on the test set, where the optimal ranks so obtained are shown in brackets. RTRMC with rank-one updating did not give better results hence, is shown omitted. The best test score of $0.7323$ is obtained by ScGrassMC at rank $10$ followed by the score $0.7370$ of {\rmc} at rank $9$. However, {\rmc} is \emph{twice} as fast as ScGrassMC.

\textbf{Remarks.} The case studies (a) to (f) are \emph{challenging} instances of (\ref{eq:matrix_completion}) as they combine ill-conditioning and low sampling in the data. Even though these studies are not fully exhaustive, they show a general trend of the performance of different algorithms. The conclusions drawn from each case study are based on a number of runs. Each figure, however, shows a \emph{typical} instance. Similarly, even though convergence of the algorithms is shown to high accuracies, the conclusions drawn are equally valid for smaller accuracies. {\rmc} has shown faster and better convergence in all the examples. 

\section{Conclusion}
We have presented {\rmc}, an efficient Riemannian conjugate-gradient algorithm for the low-rank matrix completion problem. The algorithm stems from a novel Riemannian quotient geometry endowed with a tailored Riemannian metric on the set of fixed-rank matrices. Various numerical comparisons suggest a competitive performance of {\rmc}. At the conceptual level, the paper shows that the Riemannian optimization framework can take the advantage not only of the quotient structure of the search space of fixed-rank matrices, but also of the quadratic nature of the cost function. This viewpoint is further exploited in \cite{mishra14b}.



\section{Acknowledgment}
We thank Paul Van Dooren, Bart Vandereycken, Nicolas Boumal, and Mariya Ishteva for useful discussions.

\newpage

\bibliography{LRMC14}

\begin{thebibliography}{10}
\providecommand{\url}[1]{#1}
\csname url@samestyle\endcsname
\providecommand{\newblock}{\relax}
\providecommand{\bibinfo}[2]{#2}
\providecommand{\BIBentrySTDinterwordspacing}{\spaceskip=0pt\relax}
\providecommand{\BIBentryALTinterwordstretchfactor}{4}
\providecommand{\BIBentryALTinterwordspacing}{\spaceskip=\fontdimen2\font plus
\BIBentryALTinterwordstretchfactor\fontdimen3\font minus
  \fontdimen4\font\relax}
\providecommand{\BIBforeignlanguage}[2]{{%
\expandafter\ifx\csname l@#1\endcsname\relax
\typeout{** WARNING: IEEEtran.bst: No hyphenation pattern has been}%
\typeout{** loaded for the language `#1'. Using the pattern for}%
\typeout{** the default language instead.}%
\else
\language=\csname l@#1\endcsname
\fi
#2}}
\providecommand{\BIBdecl}{\relax}
\BIBdecl

\bibitem{markovsky08a}
I.~Markovsky, ``Structured low-rank approximation and its applications,''
  \emph{Automatica}, vol.~44, no.~4, pp. 891--909, 2008.

\bibitem{benner13a}
P.~Benner and J.~Saak, ``Numerical solution of large and sparse continuous time
  algebraic matrix {R}iccati and {L}yapunov equations: A state of the art
  survey,'' MPI Magdeburg, Tech. Rep., 2013.

\bibitem{joulin10a}
A.~Joulin, F.~Bach, and J.~Ponce, ``Discriminative clustering for image
  co-segmentation,'' in \emph{IEEE Computer Society Conference on Computer
  Vision and Pattern Recognition (CVPR)}, 2010, pp. 1943--1950.

\bibitem{rennie05a}
J.~Rennie and N.~Srebro, ``Fast maximum margin matrix factorization for
  collaborative prediction,'' in \emph{International Conference on Machine
  learning (ICML)}, 2005, pp. 713--719.

\bibitem{meyer11a}
G.~Meyer, S.~Bonnabel, and R.~Sepulchre, ``{Linear regression under fixed-rank
  constraints: a Riemannian approach},'' in \emph{International Conference on
  Machine Learning (ICML)}, 2011, pp. 545--552.

\bibitem{candes09b}
E.~J. Cand{\`e}s and B.~Recht, ``Exact matrix completion via convex
  optimization,'' \emph{Foundations of Computational Mathematics}, vol.~9,
  no.~6, pp. 717--772, 2009.

\bibitem{keshavan10a}
R.~H. Keshavan, A.~Montanari, and S.~Oh, ``Matrix completion from noisy
  entries,'' \emph{Journal of Machine Learning Research}, vol.~11, no. Jul, pp.
  2057--2078, 2010.

\bibitem{cai10a}
J.~F. Cai, E.~J. Cand\`es, and Z.~Shen, ``A singular value thresholding
  algorithm for matrix completion,'' \emph{SIAM Journal on Optimization},
  vol.~20, no.~4, pp. 1956--1982, 2010.

\bibitem{jain10a}
P.~Jain, R.~Meka, and I.~Dhillon, ``Guaranteed rank minimization via singular
  value projection,'' in \emph{Advances in Neural Information Processing
  Systems 23 (NIPS)}, 2010, pp. 937--945.

\bibitem{boumal11a}
N.~Boumal and P.-A. Absil, ``{RTRMC}: A {R}iemannian trust-region method for
  low-rank matrix completion,'' in \emph{Advances in Neural Information
  Processing Systems 24 ({NIPS})}, 2011, pp. 406--414.

\bibitem{mishra13a}
B.~Mishra, G.~Meyer, F.~Bach, and R.~Sepulchre, ``Low-rank optimization with
  trace norm penalty,'' \emph{SIAM Journal on Optimization SIAM Journal on
  Optimization}, vol.~23, no.~4, pp. 2124--2149, 2013.

\bibitem{mishra12a}
B.~Mishra, K.~Adithya~Apuroop, and R.~Sepulchre, ``A {Ri}emannian geometry for
  low-rank matrix completion,'' arXiv:1211.1550, Tech. Rep., 2012.

\bibitem{mishra14a}
B.~Mishra, G.~Meyer, S.~Bonnabel, and R.~Sepulchre, ``Fixed-rank matrix
  factorizations and {R}iemannian low-rank optimization,'' \emph{Computational
  Statistics}, vol.~29, no. 3--4, pp. 591--621, 2014.

\bibitem{ngo12a}
T.~T. Ngo and Y.~Saad, ``{Scaled gradients on Grassmann manifolds for matrix
  completion},'' in \emph{Advances in Neural Information Processing Systems 25
  (NIPS)}, 2012, pp. 1421--1429.

\bibitem{vandereycken13a}
B.~Vandereycken, ``Low-rank matrix completion by {R}iemannian optimization,''
  \emph{SIAM Journal on Optimization}, vol.~23, no.~2, pp. 1214--1236, 2013.

\bibitem{absil08a}
P.-A. Absil, R.~Mahony, and R.~Sepulchre, \emph{Optimization Algorithms on
  Matrix Manifolds}.\hskip 1em plus 0.5em minus 0.4em\relax Princeton, NJ:
  Princeton University Press, 2008.

\bibitem{nocedal06a}
J.~Nocedal and S.~J. Wright, \emph{Numerical Optimization, Second
  Edition}.\hskip 1em plus 0.5em minus 0.4em\relax New York, USA: Springer,
  2006.

\bibitem{manton02a}
J.~H. Manton, ``Optimization algorithms exploiting unitary constraints,''
  \emph{IEEE {T}ransactions on {S}ignal {P}rocessing}, vol.~20, no.~3, pp.
  635--650, 2002.

\bibitem{wen12a}
Z.~Wen, W.~Yin, and Y.~Zhang, ``Solving a low-rank factorization model for
  matrix completion by a nonlinear successive over-relaxation algorithm,''
  \emph{Mathematical Programming Computation}, vol.~4, no.~4, pp. 333--361,
  2012.

\bibitem{boumal14a}
N.~Boumal, B.~Mishra, P.-A. Absil, and R.~Sepulchre, ``Manopt: a {M}atlab
  toolbox for optimization on manifolds,'' \emph{Journal of Machine Learning
  Research}, vol.~15, no. Apr, pp. 1455--1459, 2014.

\bibitem{buchanan05a}
A.~Buchanan and A.~Fitzgibbon, ``Damped {N}ewton algorithms for matrix
  factorization with missing data,'' in \emph{IEEE Computer Society Conference
  on Computer Vision and Pattern Recognition (CVPR)}, 2005, pp. 316--322 vol.
  2.

\bibitem{ring12a}
W.~Ring and B.~Wirth, ``Optimization methods on {Riemannian} manifolds and
  their application to shape space,'' \emph{SIAM Journal on Optimization},
  vol.~22, no.~2, pp. 596--627, 2012.

\bibitem{burer03a}
S.~Burer and R.~Monteiro, ``A nonlinear programming algorithm for solving
  semidefinite programs via low-rank factorization,'' \emph{Mathematical
  Programming}, vol.~95, no.~2, pp. 329--357, 2003.

\bibitem{balzano10a}
L.~Balzano, R.~Nowak, and B.~Recht, ``Online identification and tracking of
  subspaces from highly incomplete information,'' in \emph{Annual Allerton
  Conference on Communication, Control, and Computing,}, 2010, pp. 704--711.

\bibitem{mishra14b}
B.~Mishra and R.~Sepulchre, ``{R}iemannian preconditioning,'' arXiv:1405.6055,
  Tech. Rep., 2014.

\end{thebibliography}
\bibliographystyle{IEEEtran}

\end{document}